\documentclass[11pt]{article}
\usepackage[greek,english]{babel}
\usepackage[iso-8859-7]{inputenc}
\usepackage{amssymb}
\usepackage{amsmath}
\usepackage{amsthm}
\usepackage{latexsym}
\usepackage{amsfonts}
\usepackage{graphicx}
\usepackage{graphics}

\theoremstyle{plain}
\newtheorem{thm}{Theorem}[section]   
\newtheorem{prop}[thm]{Proposition}
\newtheorem{lem}[thm]{Lemma}

\theoremstyle{definition}
\newtheorem{rem}[thm]{Remark}
\newtheorem{exm}[thm]{Example}

\newtheorem*{Proof}{Proof}

\newcommand{\bbb}[1]{\mbox{\boldmath$#1$}}

\newcommand{\OO} {{\varOmega}}

\newcommand{\bi} {{\beta}}
\newcommand{\ga} {{\gamma}}

\newcommand{\ld} {{\ldots}}
\newcommand{\sm} {{\smallsetminus}}
\newcommand{\thi} {{\theta}}

\newcommand{\de} {{\delta}}
\newcommand{\De} {{\varDelta}}

\newcommand{\Si} {{\varSigma}}
\newcommand{\la} {{\lambda}}
\newcommand{\el} {{\ell}}

\newcommand{\vPi} {{\varPi}}

\newcommand{\e} {{\varepsilon}}

\newcommand{\dis}{\displaystyle}

\newcommand{\cu}{{\cal{U}}}
\newcommand{\cp}{{\cal{P}}}
\newcommand{\ch}{{\cal{H}}}
\newcommand{\cm}{{\cal{M}}}

\newcommand{\cde}{{\cal{D}}}

\newcommand{\mfd}{{\mathfrak{D}}}
\newcommand{\ra}{{\rightarrow}}
\newcommand{\oD}{{\overline{D}}}
\newcommand{\fa}{{\forall}}
\newcommand{\ko}{{\overset{\circ}{K}}}

\newcommand{\tint}{{\text{Int}}}
\newcommand{\tind}{{\text{Ind}}}
\newcommand{\qb}{$\quad\blacksquare$}
\def\1{\it1\hspace*{-0.150cm}{\footnotesize{I}}}

\def\R{{\mathbb{R}}}
\def\C{{\mathbb{C}}}
\def\Q{{\mathbb{Q}}}
\def\N{{\mathbb{N}}}

\def\O{{\mathbb{O}}}

\begin{document}
\title{\bf Universal Taylor series on specific compact sets}
\author{\bf N. Tsirivas}\footnotetext{\hspace{-0.5cm}}
\footnotetext{{The research project is implemented within the framework of the Action ``Supporting Postdoctoral Researchers'' of the Operational Program ``Educational and Lifelong Learning'' (Action's Beneficiary: General Secretariat for Research and Technology), and is co-financed by the European Social Fund (ESF) and the Greek State.}}
\date{}
\maketitle
\noindent
{\bf Abstract:} Let $D$ be the unit disc.\ We denote by $\C$ the set of complex numbers and consider the set $\cm_D:=\{K\in\cp(\C)\mid K$ is compact, $K^c$ is connected, $K\cap D=\emptyset\}$.\ Let $A(K):=\{f:K\ra\C\mid f$ is continuous and holomorphic in $\overset{\circ}{K}\}$, for $K\in\cm_D$.\ The space $A(K)$, for $K\in\cm_D$, is endowed with the supremum norm. It is a well known result \cite{20} that there exist holomorphic functions $f$ on $D$ for which the partial sums $S_n(f)$, $n=1,2,\ld$ of the Taylor series with center $0$ are dense in $A(K)$ for every $K\in\cm_D$.
It is also known that the above result fails \cite{23} if we consider the weighted polynomials $2^nS_n(f)$, $n=1,2,\ld$ instead of $S_n(f)$, $n=1,2,\ld\;.$ In the opposite direction, the main result of this work shows that there exist holomorphic functions $f$ on $D$ for which the sequence $2^nS_n(f)$, $n=1,2,\ld$ is dense in $A(K)$ for specific $K\in\cm_D$. In this case the geometry of $K$ plays a crucial role. We also generalize these results on arbitrary simply connected domains.
\vspace*{0.2cm} 
\vspace*{0.2cm} 
\section{Introduction}\label{sec1}
\noindent

Let $\OO\subseteq\C$ be an open set, $z_0\in\OO$ and $\ch(\OO):=\{f:\OO\ra\C\mid f\;\text{is holomorphic}\}$. For $f\in\ch(\OO)$, the symbol  $S_n(f,z_0)$ denotes the $n$-th partial sum of the Taylor's development of $f$ with center $z_0$, that is
\[
S_n(f,z_0)(z):=\sum^n_{k=0}\frac{f^{(k)}(z_0)}{k!}z^k, \ \ z\in\C, \ \ n=0,1,2,\ld\;.
\]
The open unit disk $\{z\in\C\mid|z|<1\}=D(0,1)$ is denoted by $D$.
Vassili Nestoridis proved in the very influential paper \cite{20} the following result:

\textit{there exists $f\in\ch(D(0,1))$ such that for every compact set $K\subset D^c$, with connected complement and for every $h:K\ra\C$, $h$ continuous and holomorphic in $\overset{\circ}{K}$, there exists a sequence of natural numbers $(\la_n)$ such that}
\[
 S_{\la_n}(f,0)\ra h \ \ \text{as} \ \ n\ra+\infty  \ \ \ \textit{uniformly on} \ K.
\]

The set of functions $f\in\ch(D)$ that satisfy the above property is the well known set of universal Taylor series, denoted by $\cu(D,0)$. In \cite{23} the present author examined the above problem in a more general frame. More specifically, he fixed a sequence $(\bi_n)$ of complex numbers and examined whether an approximation scheme as above may hold for the weighted partial sums $(\bi_nS_n(f,0))$. Namely, does there exist $f\in\ch(D)$ such that for every $K$ and $h$ as above,
\[
\bi_{\la_n}S_{\la_n}(f,0) \ra h \ \ \text{as} \ \ n\ra+\infty  \ \ \ \text{uniformly on} \ K
\]
for some sequence of natural numbers $(\la_n)$? A complete answer to this question is given in \cite{23}, where it is proved that: \textit{the answer is positive if and only if the sequence $(|\bi_n|^{\frac{1}{n}})$ has 1 as a limit point}.

To state our main result we need to introduce a little notation and terminology. Let $(\bi_n)$ be a given sequence of complex numbers and let $K\subset D^c$ be a compact set with connected complement. Define
\[
\cu(0, D,K, (\bi_n)):=\big\{f\in\ch(D)\mid\overline{\{\bi_nS_n(f,0),\;n=1,2,\ld\}}=A(K)\big\},
\]
where $A(K):=\{f:K\ra\C\mid f$ is continuous and holomorphic in $\overset{\circ}{K}\}$ and the closure of the partial sums is taken with respect to the supremum norm on $K$.

In view of the above, suppose that $1$ is not a limit point of the sequence $(|\bi_n|^{\frac{1}{n}})$. Now the following question arises naturally. Are there compact sets $K\subset D^c$ with connected complement so that the class $\cu(0, D,K, (\bi_n))$ is non-empty?

For a given compact set $L$ with finitely many connected components we assign a number $\rho_L$, which is characteristic for the compact set, and it is called the asymptotic convergence factor of $L$. This number is always between $0$ and $1$ and for its definition see Section 2 of the present paper.
For a sequence $(\lambda_n)$ of complex numbers we denote by $L((\lambda_n)_{n=1}^{\infty})$ the set of limit points of the sequence.
We prove the following
\begin{thm} \label{mainresult}
Fix a sequence $(\bi_n)$ of complex numbers. Let $K\subset \overline{D}^c$ be compact having the following properties: $K$ has more than one elements, $K$ is connected and $K^c$ is connected. Set $L:=\overline{D}\cup K$ and consider the positive number $M:=e^{\max_{z\in \overline{D}}g_{\Omega}(z,\infty) }$, where $\Omega :=(\mathbb{C}\cup \{ \infty \} )\setminus K$ and $g_{\Omega }$ is the Green function for $\Omega $ with pole at infinity.\\
\begin{enumerate}
\item[(i) ] If $L((|\bi_n|^{\frac{1}{n}})_{n=1}^{\infty}) \cap (\rho_L, \frac{1}{\rho_L} ) \neq \emptyset$  then $\cu(0, D,K, (\bi_n))\neq \emptyset$.
\item[(ii)] If $\limsup_n|\bi_n|^{\frac{1}{n}}<\frac{1}{\textrm{dist}(0,K)} $
then $\cu(0, D,K, (\bi_n))=\emptyset$.
\item[(iii)] If $\liminf_n|\bi_n|^{\frac{1}{n}}>M$ then $\cu(0, D,K, (\bi_n))=\emptyset$.
\end{enumerate}
\end{thm}

Some remarks are in order.\\
Firstly, observe that the numbers $\rho_L$, $\textrm{dist}(0,K)$, $M$ involved in Theorem \ref{mainresult} depend only on the sets $K,D$. It is rather surprising that as a consequence of Theorem \ref{mainresult}, which concerns the occurrence of universality of the partial sums on certain compact sets $K$, we extract information on the relation between the above numbers. In particular, it follows that
$$ \frac{1}{\textrm{dist}(0,K)} <\rho_L <\frac{1}{\rho_L} <M .$$
These inequalities can also be proved by the aid of potential theory, without the use of universality; however, this approach is more involved. For details see \cite{24}.\\
Secondly, if one wishes to know all the complex sequences $(\bi_n)$ and all the appropriate compact sets for which the set $\cu(0, D,K, (\bi_n))$ is non-empty, then one has to deal with the cases which are not covered by Theorem \ref{mainresult}. For instance, the following question remains open.\\
{\bf{Question:}}
\textit{is true that $\cu(0, D,K, (\bi_n))=\emptyset$ provided that}
$$ L((|\bi_n|^{\frac{1}{n}})_{n=1}^{\infty}) \subset \left[\frac{1}{\textrm{dist}(0,K)} ,\rho_L \right] \cup \left[\frac{1}{\rho_L} ,M\right] ?$$
Thirdly, we observe that the compact set $K$, appearing in Theorem \ref{mainresult}, is not a singleton. This fact is due to technical reasons. In particular, the use of potential theory forces us to consider ``fat" sets, that is, the domain $L^c$ $(L:=\overline{D}\cup K)$ should be regular, see \cite{21}. This means that if one wants to include in Theorem \ref{mainresult} sets $K$ consisting of finitely many elements then a different approach is needed.\\ 

To illustrate Theorem \ref{mainresult} we present some examples. By $D(z,r), \overline{D(z,r)}$ we denote the open, closed disk with center $z$ and radius $r>0$ respectively. Consider the sequence $\bi_n=2^n$, $n=1,2,\ldots $. Then $\cu(0, D, \overline{D(a,1)}, (2^n))\neq \emptyset$ provided that $a$ is positive real number with $a\geq 18$. On the other hand, for every closed disk $\overline{D(z,r)}$ with $\overline{D(z,r)} \subset \{ w\in \mathbb{C}: 1<|w|<2 \}$ we have $\cu(0, D, \overline{D(z,r)}, (2^n))= \emptyset$. The above facts are explained in great detail in Sections \ref{sec3},\ref{sec4}.\\   

Our paper is organized as follows. In Section \ref{sec2} we present the necessary terminology and we prove item $(i)$ of Theorem \ref{mainresult}. A sample of examples is given in Section \ref{sec3}.
In Section \ref{sec4} we prove items $(ii)$ and $(iii)$ of Theorem \ref{mainresult}. Actually, we shall prove a stronger version of Theorem \ref{mainresult}, see Theorem \ref{thm2.2}, Proposition \ref{prop4.1}, Proposition \ref{prop4.2},  which covers more general domains and not only the unit disk. The methods in this paper come from potential theory. Potential Theory has recently become a powerful tool in solving a variety of problems on universal Taylor series; see the respective papers \cite{5}, \cite{6}, \cite{8}, \cite{10} - \cite{14}, \cite{17}, \cite{18}, \cite{19}, \cite{23}, \cite{25}.

An abstract theory of universal series is developed in \cite{3}, that covers the result of V. Nestoridis in \cite{20}. Later in \cite{5}, \cite{16}, \cite{23}, the authors extended the above theory. This line of research is closely related to our investigations. For this reason our paper can be considered as a continuation of this series of papers.
\section{The main result}\label{sec2}
\noindent

First of all we develop here the necessary terminology of our paper.

We define now the asymptotic convergence factor $\rho_L$ for a compact set of the form $L=\bigcup\limits^{m}_{i=0}K_i$.

Let $L$ be a non-connected compact subset of $\C$, with connected complement.

We further assume that $L:=\bigcup\limits^m_{i=0}K_i$, for some $m\in\N$, $m\ge1$ where $K_i$, $i=0,1,\ld,m$ are the connected components of $L$. Let $p_i$, $i=0,1,\ld,m$ be $m+1$ different complex polynomials, that is $p_i\neq p_j$ for every $i,j\in\{0,1,\ld,m\}$, $i\neq j$.

We consider the function $F:L\ra\C$, that is defined by the formula:
\[
F(z)=p_j(z) \ \ \text{if} \ \ z\in K_j, \ \ \text{for every} \ \ j\in\{0,1,\ld,m\}.
\]
We fix some positive number $\de$. The problem is to find a polynomial $p$ such that $\|F-p\|<\de$ and the relation between $p$ and $\de$ in any case.

For every $n=1,2,\ld$, let $V_n$ be the set of complex polynomials with degree at most $n$.

We denote
\[
d_n, F:=\min\{\|F-p\|,\;p\in V_n\} \ \ \text{for} \ \ n=1,2,\ld.
\]
Of course for every $n\in\N$, there exists some $p\in V_n$ such that $d_n,F=\|F-p\|_L$, and the polynomial $p$ is unique (\cite{26}) for every $n\ge1$. Even, if the formulation of the problem of finding the above best polynomial $p$ that minimizes the quantity $\|F-p\|_L$ is simple this is usually unknown and difficult to compute (see \cite{9}, page 11).

However, if the compact set $L$ has a simple construction and good properties the previous approximation problem can be solved. However, the computation of the best polynomial is difficult even if for simple cases and in most of cases this is become with numerical methods that are complicated.

A classical theorem in this region is the following.
\begin{thm}\label{thm2.1}
The number $\rho_L:=\underset{n\ra+\infty}{\lim\sup}\,d^{\frac{1}{n}}_n,F$ is a positive constant such that $\rho_L\in(0,1)$ and is independent from the function $F$ and it is dependent only on the compact set $L$.
\end{thm}

The number $\rho_L$ is called the asymptotic convergence factor of $L$ and is a characteristic for the compact set $L$.

For the topological concepts of this paper we refer to the classical book of Burckel \cite{4}. More specific for the definitions of a curve, or a loop, or an arc, or a simple curve, or a smooth curve see Definition 1.11 \cite{4}.

With a Jordan curve we mean a homeomorphism in $\C$ of a circle.

If $\ga$ is a smooth Jordan curve and $w\in\C\backslash\ga$, the index $\tind_\ga(w):=\dfrac{1}{2\pi i}\dis\int\limits_\ga\dfrac{1}{z-w}dz$. For a compact subset $K$ of $\C$ and a Jordan curve $\ga$ such that $\ga\cap K=\emptyset$, we write
\[
\tind_\ga(K):=\{\tind_\ga(w),\;w\in K\}.
\]
The definition of interior, $\tint(\ga)$, and Exterior $Ex(\ga)$ of a Jordan curve $\ga$ is given in Definition 4.45 (i) of \cite{4}. For results about potential theory we refer to the classical books \cite{1} and \cite{21}.

We consider a compact set $L\subseteq\C$, with $L^c$ connected such that $L=\bigcup\limits^{m_0}_{i=0}K_i$, $m_0\in\N$, where $K_i$, $i=0,1,\ld,m_0$ be the connected components of $L$.

We consider the set:

$\mfd_L:=\{\De\in\cp(\C)\mid$ there exist $m_0+1$ smooth Jordan curves $\de_i$, for $i=0,1,\ld,m_0$ such that $\De=\bigcup\limits^{m_0}_{i=0}\de_i$, $K_i\subset\tint(\de_i)$ and $\tind\de_i(K_i)=\{1\}$ for every $i=0,1,\ld,m_0$ and $\bigcup\limits^{m_0}_{i=0\atop i\neq j}\de_i\subset Ex(\de_j)$ for every $j=0,1,\ld,m_0\}$.

By Lemma 1.2 of \cite{24} we have the $\mfd_L\neq\emptyset$ and the set is uncountable. From now on and till the end of the paper we consider a compact set $L=\bigcup\limits^{m_0}_{i=0}K_i$, $m_0\in\N$, where $K_i$, $i=0,1,\ld,m_0$ are the connected components of $L$, $L^c$ is connected, $\overset{\circ}{K}_0\neq\emptyset$, and the compact sets $K_i$, $i=1,\ld,m_0$ contain more than one point.

Let $\OO:=(\C\backslash L)\cup\{\infty\}$. Then there exists the unique Green's function $g_\OO$ for $\OO$, with pole at infinity (Definition 4.4.1 and Theorem 4.4.2 of \cite{21}).

Let $\De\in\mfd_L$. We write
\[
\thi_{L,\De}:=\max e^{-g_\OO(z,\infty)}:=\max\{x\in\R\mid\exists\;z\in\De:x=e^{-g_\OO(z,\infty)}\}
\]
and
\[
\thi_L:=\inf\{x\in\R\mid\exists\;\De\in\mfd_L:x=\thi_{L,\De}\}.
\]
In \cite{24} is proved that $\thi_L\in(0,1)$ and by Proposition 2.3 of \cite{24} we have that $\rho_L=\thi_L$.

We define $\vPi:=L\sm K_0=\bigcup\limits^{m_0}_{i=1}K_i$.
We fix some sequence $(\bi_n)$ of complex numbers and we also fix some point $z_0\in\ko_0$.
Let $f\in A(K_0)$. For $z\in\C$, $S_n(f,z_0)(z)$ denotes the $n$-th partial sum of the Taylor's development of $f$ with center $z_0$, for $n=0,1,2,\ld\;$.
Now the set of universal vectors of $A(K_0)$ with respect to the point $z_0$, the compact set $\vPi$ and the sequence $(\bi_n)$ is defined to be the set
\[
\cu(z_0,K_0,\vPi,(\bi_n)):=\cu(\bi_n)=\big\{
f\in A(K_0)\mid\overline{ \{ \bi_nS_n(f,z_0) :n=1,2,\dots \} }=A(\vPi)\big\}.
\]
The main problem of this paper is to find sequences $(\bi_n)$ such that $\cu(\bi_n)\neq\emptyset$. We write only the dependence of the sequence $(\bi_n)$ in the symbol $\cu(\bi_n)$ because all the others parameters of this are supposed to be fixed, that is the point $z_0$, the compact sets $K_0$, $\vPi$ and of course the partial sums $S_n(f,z_0)$, for some $f\in A(K_0)$, $n=0,1,2,\ld\;.$ We will see that there are many sequences $(\bi_n)$ such that $\cu(\bi_n)\neq\emptyset$, in general.

Firstly, let us examine the candidate sequences $(\bi_n)$.

We remark that if the sequence $(\bi_n)$ is finally zero, that is if there exists some natural number $n_0\in\N$ such that $\bi_n=0$ for every $n\ge n_0$ then $\cu(\bi_n)=\emptyset$. So we have to examine the case where there exists an infinity number of integers $n\in\N$ such that $\bi_n\neq0$.

Firstly we examine the case where $\bi_n\neq0$ for every $n=0,1,2,\ld\;.$

Now we examine two cases: \vspace*{0.2cm} \\
{\bf First case:} The sequence $(\bi_n)$ has a finite limit point non-zero. This means that there exists some subsequence $(\bi_{k_n})$ of $(\bi_n)$ and some  complex number $w_0\neq0$ such that $\bi_{k_n}\ra w_0$ as $n\ra+\infty$.

Then $|\bi_{k_n}|^{1/k_n}\ra1$ and by the main result of \cite{23} we have of course $\cu(\bi_n)\neq\emptyset$. \vspace*{0.3cm}
{\bf Second case:} The sequence $(\bi_n)$ does not have a finite non-zero limit point. Of course the sequence $(\bi_n)$ as a sequence of complex numbers has obligatory limit points on $\C\cup\{\infty\}$. So the sequence $(\bi_n)$ has only two possible limit points: 0 and $\infty$.

We consider firstly the case where $\bi_n\ra0$. In this case we have $\underset{n\ra+\infty}{\lim\sup}|\bi_n|^{1/n}\le1$.

If $\underset{n\ra+\infty}{\lim\sup}|\bi_n|^{1/n}=1$ then by the main result of \cite{23} we have that $\cu(\bi_n)\neq\emptyset$. So the interesting case is to examine the case where $\underset{n\ra+\infty}{\lim\sup}|\bi_n|^{1/n}\in[0,1)$.

We consider now the case where $\bi_n\ra\infty$. Then it is obvious that
$\underset{n\ra+\infty}{\lim\inf}|\bi_n|^{1/n}\ge1$. If $\underset{n\ra+\infty}{\lim\inf}|\bi_n|^{1/n}=1$ then by the main result of \cite{23} we have that $\cu(\bi_n)\neq\emptyset$. So the interesting case here is when $\underset{n\ra+\infty}{\lim\inf}|\bi_n|^{1/n}\in(1,+\infty]$.

The case where the sequence $(\bi_n)$ has exactly two limit points 0 and $\infty$ is reduced on the two above cases passing to a suitable subsequence of $(\bi_n)$.

So, in all the above cases the problem is reduced in two only cases:

1) $\bi_n\ra0$ and $|\bi_n|^{1/n}\ra a_0\in[0,1)$

2) $\bi_n\ra\infty$ and $|\bi_n|^{1/n}\ra a_0\in(1,+\infty]$

We will prove the following result, which establishes item $(i)$ of Theorem \ref{mainresult} 
\begin{thm}\label{thm2.2}
By the previous terminology we consider a sequence $(\bi_n)$ of complex numbers such that the sequence $|\bi_n|^{1/n}$ has a limit point in the open interval $\Big(\rho_L,\dfrac{1}{\rho_L}\Big)$. Then the set $\cu(\bi_n)$ is a dense $G_\de$-subset of $A(K_0)$. As a consequence the set $\cu(\bi_n)\neq\emptyset$.
\end{thm}

We will prove Theorem \ref{thm2.2} using some lemmas and Baire's Category Theorem.
For this we write $\cu(\bi_n)=\bigcap\limits^{+\infty}_{n=1}V_n$ where each $V_n$ is open and dense subset of $A(K_0)$. Let us describe this procedure as follows.
Let $(u_j)$, $j=1,2,\ld$ be an enumeration of all non-zero complex polynomials with coefficients in $\Q+i\Q$. For every $j,s,n\in\N$ define:
\[
V(j,s,n):=\{f\in A(K_0)\mid\|\bi_nS_n(f,z_0)-u_j\|_\vPi<1/s\},
\]
where
\[
\|h\|_\vPi:=\max\{|h(z)|:z\in K\}\}.
\]
for every $h\in A(\vPi)$.
\begin{lem}\label{lem2.3}
The sets $V(j,s,n)$ are open subsets of $(A(K_0),\|\cdot\|_\infty)$ for every $j,s,n\in\N$.
\end{lem}
\begin{lem}\label{lem2.4}
The following holds:
\[
\cu(\bi_n)=\bigcap^{+\infty}_{j=1}\bigcap^{+\infty}_{s=1}\bigcup^{+\infty}_{n=1}
V(j,s,n).
\]
\end{lem}

For the proofs of the above two lemmas see Lemmas 2.2 and 2.3 of \cite{23} or Lemma \ref{lem2.4} and Proposition \ref{prop2.5} of \cite{16}, and using Mergelyan's Theorem \cite{22}. We complete the proof of Theorem \ref{thm2.2} with the following proposition.
\begin{prop}\label{prop2.5}
By the suppositions of Theorem \ref{thm2.2} the sets $\bigcup\limits^{+\infty}_{n=1}V(j,s,n)$ are dense in $(A(K_0),\|\cdot\|_\infty)$ for every $j,s\in\N$.
\end{prop}
\begin{Proof}
It suffices to examine the case where $\bi_n\ra0$ and $\rho_L<\dis\lim_{n\ra+\infty}|\bi_n|^{1/n}<1$ where $(\bi_n)$ is a sequence of non-zero complex numbers.

We fix some natural numbers $j_0$ and $s_0$ and we write $B:=\bigcup\limits^{+\infty}_{n=1}V(j_0,s_0,n)$. We show that the set $B$ is dense in $(A(K_0),\|\cdot\|_\infty)$.

By Mergelyan's Theorem (see \cite{22}) the set of all complex polynomials is dense in $A(K_0)$.

So, it suffices to show that the set $B$ has non-empty intersection with every open neighbourhood of every polynomial.

So we fix some polynomial $p_0$ and some positive number $\e_0>0$.

We consider the open neighbourhood of $p_0$ in the space $(A(K_0),\|\cdot\|_\infty)$
\[
V_{p_0,\e_0}:=\{h\in A(K_0)\mid\|p_0-h\|_{K_0}<\e_0\}.
\]
It suffices to show that $B\cap V_{p_0,\e_0}\neq\emptyset$.

This means that there exists some complex function $f\in A(K_0)$ and some natural number $N_0$ such that:
\setcounter{equation}{0}
\begin{eqnarray}
\|f-p_0\|_{K_0}<\e_0 \ \ \text{and} \label{sec31}
\end{eqnarray}
\begin{eqnarray}
\|\bi_{N_0}S_{N_0}(f,z_0)-u_{j_0}\|_\vPi<\frac{1}{s_0}.  \label{sec32}
\end{eqnarray}
By our hypothesis for the sequence $(\bi_n)$ we have:
\[
\rho_L<\lim_{n\ra+\infty}|\bi_n|^{1/n}<1.
\]
We write $b_0:=\dis\lim_{n\ra+\infty}|\bi_n|^{1/n}$. Then $b_0\in(\rho_L,1)$.
This means that there exists some strictly increasing subsequence $(\la_n)$ of natural numbers such that:
\[
|\bi_{\la_n}|^{1/\la_n}\ra b_0 \ \ \text{as} \ \ n\ra+\infty.
\]
So there exists some positive number $\e_1$ and some natural number $v_1\in\N$ such that:
\begin{eqnarray}
\rho_L<b_0-\e_1<|\bi_{\la_n}|^{1/\la_n}<b_0+\e_1<1  \label{sec33}
\end{eqnarray}
for every $n\in\N$, $n\ge v_1>2$.

We choose some positive number $c_0\in(\rho_L,b_0-\e_1)$ and we fix it. After the choice of the positive number $c_0$ we choose some $\De_0\in\mfd_L$ that depends on $c_0,L$ such that
\[
\thi_L=\rho_L<\thi_{L,\De_0}<c_0
\]
by the definition of $\thi_L$ and Proposition 2.6 of \cite{24}.

Now, for every $n\ge\la_{v_1}$ we consider the well defined complex function $F_n:L\ra\C$ as follows:
\[
F_n(z):=\left.\begin{array}{lcc}
                p_0(z) & \text{if} & z\in K_0 \\ [1.5ex]
                \dfrac{1}{\bi_n}u_{j_0}(z) & \text{if} & z\in\vPi
              \end{array}\right\}
\]
For every function $F_n$, $n\ge\la_{v_1}$ we apply Proposition 2.2 of \cite{24} and we take that there exists some natural number $m_0=m_{L,\De_0,c_0}$ independent from $n$, some positive constant $A_n=A_{L,\De_0,c_0,F_n}$ and some sequence $S^n_m$ of polynomials that depends on $L,\De_0,c_0,F_n$ such that:
\begin{eqnarray}
\|F_n-S^n_m\|_L<A_n\cdot c^m_0 \ \ \text{for every} \ \ m\in\N, \ \ m\ge m_{L,\De_0,c_0}, \ \ n\ge\la_{v_1}.  \label{sec34}
\end{eqnarray}
By (\ref{eq4}) for every $n=m=\la_v$, where $v>\max\{m_{L,\De_0,c_0},v_1\}=v_2$ we get
\begin{eqnarray}
\|F_{\la_v}-S^{\la_v}_{\la_v}\|_L<A_{\la_v}\cdot c^{\la_v}_0.  \label{sec35}
\end{eqnarray}
By the proof of Proposition 2.2 \cite{24} we have:
\[
A_{\la_v}:=\frac{\la_0\cdot\|F_{\la_v}\|_L}{2\pi\cdot dist(\De_0,L)},
\]
where the number $C_1:=\dfrac{\la_0}{2\pi dist(\De_0,L)}$ is a positive constant that depends only on $L,\De_0,c_0$ and it is independent from $v$.

So we have
\[
A_{\la_v}=C_1\cdot\|F_{\la_v}\|_L \ \ \text{for every} \ \ v\ge v_2.
\]
We have $u_{j_0}\neq0$. So there exists some natural number $v_3$ such that $\dfrac{\|u_{j_0}\|_\vPi}{|\bi_n|}>\|p_0\|_{K_0}$ for every $n\ge v_3$ because $\bi_n\ra0$ as $n\ra+\infty$. By the above we have:
\begin{eqnarray}
A_{\la_v}=C_1\cdot\frac{\|u_{j_0}\|_\vPi}{|\bi_{\la_v}|} \ \ \text{for every} \ \ n\ge v_4:=\max\{v_3,v_2\}.  \label{sec36}
\end{eqnarray}
By (\ref{sec35}) and (\ref{sec36}) we get:
\begin{eqnarray}
\|F_{\la_v}-S^{\la_v}_{\la_v}\|_L<C_1\cdot\frac{\|u_{j_0}\|_\vPi}{|\bi_{\la_v}|}\cdot
c^{\la_v}_0 \ \ \text{for every} \ \ v\ge v_4.  \label{sec37}
\end{eqnarray}
We apply (\ref{eq7}) on the compact sets $K_0$ and $\vPi$ and we take:
\begin{eqnarray}
\|S^{\la_v}_{\la_v}-p_0\|_{K_0}<C_1\cdot\frac{\|u_{j_0}\|_\vPi}{|\bi_{\la_v}|}\cdot
c^{\la_v}_0 \ \ \text{for} \ \ v\ge v_4  \label{sec38}
\end{eqnarray}
and
\begin{eqnarray}
\bigg\|S^{\la_v}_{\la_v}-\frac{1}{\bi_{\la_v}}u_{j_0}\bigg\|_\vPi<C_1\cdot
\frac{\|u_{j_0}\|_\vPi}{|\bi_{\la_v}|}\cdot c^{\la_v}_0 \ \ \text{for} \ \ v\ge v_4.  \label{sec39}
\end{eqnarray}
By (\ref{sec39}) we get:
\begin{eqnarray}
\|\bi_{\la_v}\cdot S^{\la_v}_{\la_v}-u_{j_0}\|_\vPi<C_1\cdot\|u_{j_0}\|_\vPi\cdot c^{\la_v}_0, \ \ v\ge v_4.  \label{sec310}
\end{eqnarray}
By (\ref{sec33}) we have:
\begin{eqnarray}
\frac{1}{|\bi_{\la_v}|}<\frac{1}{(b_0-\e_1)^{\la_v}} \ \ \text{for} \ \ v\ge v_1.  \label{sec311}
\end{eqnarray}
Because $v_4\ge v_1$ by (\ref{eq8}) and (\ref{sec311}) we have:
\begin{eqnarray}
\|S^{\la_v}_{\la_v}-p_0\|_{K_0}<C_1\|u_{j_0}\|_\vPi\cdot\bigg(\frac{c_0}{b_0-\e_1}\bigg)^{\la_v}, \ \ v\ge v_4  \label{sec312}
\end{eqnarray}
we have $\dfrac{c_0}{b_0-\e_1}\in(0,1)$ because the choice of the positive number $c_0$.

By (\ref{sec310}) and (\ref{sec312}) if we take $v_0\ge v_4$ big enough we can take that:
\begin{eqnarray}
\|S^{\la_{v_0}}_{\la_{v_0}}-p_0\|_{K_0}<\e_0 \ \ \text{and} \label{sec313}
\end{eqnarray}
\begin{eqnarray}
\|\bi_{\la_v}S^{\la_{v_0}}_{\la_{v_0}}-u_{j_0}\|_\vPi<\frac{1}{s_0}. \label{sec314}
\end{eqnarray}
The polynomial $S^{\la_{v_0}}_{\la_{v_0}}$ is a polynomial of degree at most $\la_{v_0}-1$. This gives that $S_{\la_{v_0}}(S^{\la_{v_0}}_{\la_{v_0}},z_0)=S^{\la_{v_0}}_{\la_{v_0}}$ and setting $f:=S^{\la_{v_0}}_{\la_{v_0}}$ we have satisfied (\ref{eq1}) and (\ref{eq2}) for $N_0:=\la_{v_0}$, because $f=S^{\la_{v_0}}_{\la_{v_0}}\in A(K_0)$ as a polynomial.

So we have proved completely the case where $\bi_n\ra0$ and $\rho_L<\dis\lim_{n\ra+\infty}|\bi_n|^{1/n}<1$.

The case where $\bi_n\ra\infty$ and $\underset{n\ra+\infty}{\lim\inf}|\bi_n|^{1/n}\in\Big(1,\dfrac{1}{\rho_L}\Big)$ is very similar with the above case and for this its proof is omitted and is left as an easy exercise for the interested reader. In the case where the sequence $\bi_n$ has exactly two limit points, 0 and $\infty$ we apply almost the same the above proof for two suitable subsequences of $(\bi_n)$ that tend to 0 and $\infty$.

Finally, in the case where the sequence $(\bi_n)$ has infinite terms $\bi_n$ such that $\bi_n=0$, we pass to a suitable subsequence of $(\bi_n)$ with non-zero terms and we apply again the above proof for this subsequence. This completes the proof of this proposition. \qb
\end{Proof}
\noindent
{\bf Proof of Theorem \ref{thm2.2}.} By Lemmas \ref{lem2.3}, \ref{lem2.4}, Proposition \ref{prop2.5}, the fact that the space $(A(K_0),\|\cdot\|_\infty)$ is a complete metric space and Baire's Category Theorem we conclude Theorem \ref{thm2.2}. \qb\vspace*{0.2cm}

By Theorem \ref{thm2.2} we have found that there exist many sequences $(\bi_n)$ such that $\cu(\bi_n)\neq\emptyset$.

However, the number $\rho_L$ even if is a very well known, in the literature,constant for $L$, it is very hard to be computed in many cases. This means that the value of Theorem \ref{thm2.2} is existential in part.

Instead, we give specific examples of sequences $(\bi_n)$ such that $\cu(\bi_n)\neq\emptyset$ in the following Section \ref{sec3}.
\section{Some specific examples of sequences such that $\bbb{\cu(\bi_n)\neq\emptyset}$}\label{sec3}
\begin{exm}\label{exm3.1}
Even if it is very hard to compute with some accuracy the number $\rho_L$, in general for a compact set $L$, Theorem \ref{thm2.2} is strong enough in order to give us specific examples of sequences $(\bi_n)$ such that $\cu(\bi_n)\neq\emptyset$.

Of course we search sequences $(\bi_n)$ such that the number 1 is not a limit point of the sequence $\Big(\sqrt[n]{|\bi_n|}\Big)$, because the other cases are analysed completely in \cite{23}. So the problem is the following: To give specific examples of compact sets $L$ where $L=\bigcup\limits_{i=0}K_i$, with the above properties as in Section \ref{sec2} and sequences $(\bi_n)$ of complex numbers such that
\[
\cu(z_0,K_0,\vPi,(\bi_n))=\cu(\bi_n)\neq\emptyset, \ \ \vPi:=L\backslash K_0.
\]
In order to give an example using Theorem \ref{thm2.2} it is not need to compute the number $\rho_L$ for some ``good'' compact set $L$ but to find an arbitrary number $\la\in(\rho_L,1)$, because then the sequence $\bi_n=\la^n$, $n=1,2,\ld$ gives us an example. So we have two choices. If we fix the compact set $L$ to find an upper bound of the number $\rho_L$, $\la$, where $\la<1$, or if we fix a number $\la\in(0,1)$ to find a compact set $L$ such that $\rho_L<\la$. The second choice is much more easier than the first, so we apply it.

For example, let us find an example of some compact set $L$ such that the respective sequence to be the sequence $\bi_n=2^n$, $n=1,2,\ld\;$.

By Theorem \ref{thm2.2} it suffices to find a compact set $L$ such that $\rho_L<\dfrac{1}{2}$. The simplest case of course is the case where the set $L$ is consisted from two disjoint simple compact sets $K_0$ and $K_1$. For simplicity we take $K_0:=\oD(0,1)$, where $\oD(0,1):=\{z\in\C\mid\,\mid z\mid\le1\}$ and
\[
K_1:=\oD(\thi_0,1):=\{z\in\C\mid\,\mid z-\thi_0\mid\le1\} \ \ \text{for some} \ \ \thi_0>1.
\]
So, the problem is if there exists $L=K_0\cup K_1$, where $K_0,K_1$ are as above, $K_0\cap K_1=\emptyset$ and $\rho_L<\dfrac{1}{2}$.

We will see that this can happen.

We consider the polynomial $p$ where $p(z):=z(z-\thi_0)$.

For $z\in K_0$ we have $|z|\le1$ and $|z-\thi_0|\le1+\thi_0$, so $\|p\|_{K_0}\le1+\thi_0$. But we have $|p(-1)|=1+\thi_0$. So $\|p\|_{K_0}=1+\thi_0$.

Similarly, we take that $\|p\|_{K_1}=1+\thi_0$. Thus, we have $\|p\|_L=1+\thi_0$.

Let
\[
C_1:=C\bigg(0,\frac{\thi_0}{2}-1\bigg):=\bigg\{z\in\C\mid\,\mid z\mid=\frac{\thi_0}{2}-1\bigg\}.
\]
We want to have: $\dfrac{\thi_0}{2}-1>1\Leftrightarrow\thi_0>4$. So, we suppose that $\thi_0>4$. This gives that the circle $C_1$ has the set $K_0$ in its interior.

Let some $z\in C_1$. Then there exists $t\in[0,1]$ such that $z=\Big(\dfrac{\thi_0}{2}-1\Big)\cdot e^{2\pi it}$. We take that
\[
|p(z)|\ge\bigg(\frac{\thi_0}{2}\bigg)^2-1>0, \ \ \text{for} \ \ z\in C_1
\]
and for $t=0$ we take that
\[
|p(z_0)|=\bigg(\frac{\thi_0}{2}\bigg)^2-1 \ \ \text{for} \ \ z_0=\frac{\thi_0}{2}-1\in C_1.
\]
Let $C_2:=C\Big(\thi_0,\dfrac{\thi_0}{2}-1\Big)$. We take that the circle $C_2$ has the set $K_1$ in its interior and the circles $C_1$ and $C_2$ have each other in its exterior, and of course the circles $C_1$ and $C_2$ are smooth Jordan curves.

We set $\De:=C_1\cup C_2$, and of course $\De\in\cde_L$, by the above.

As in $C_1$ we have that
\[
|p(z)|\ge\bigg(\frac{\thi_0}{2}\bigg)^2-1 \ \ \text{for every} \ \ z\in C_2 \ \ \text{and} \ \ \min_{z\in C_2}|p(z)|=\bigg(\frac{\thi_0}{2}\bigg)^2-1.
\]
So
\[
\min_{z\in\De}|p(z)|=\bigg(\frac{\thi_2}{2}\bigg)^2-1.
\]
We set $\OO:=(\C\cup\{\infty\})\setminus L$. Let $g_\OO$ to be the Green's function for $\OO$ with pole at infinity.

By Bernstein's Lemma, Theorem 5.5.7 of \cite{21} we have that
\[
\bigg(\frac{|p(z)|}{\|p\|_L}\bigg)^{1/2}\le e^{g_\OO(z,\infty)} \ \ \text{for} \ \ z\in\OO\setminus\{\infty\}.
\]
So, because $\De\subset\OO\setminus\{\infty\}$ we take that
\[
e^{-g_\OO(z,\infty)}\le\bigg(\frac{\|p\|_L}{|p(z)|}\bigg)^{1/2} \ \ \text{for every} \ \ z\in\De.
\]
This gives that
\begin{align*}
\sup_{z\in\De}e^{-g_\OO(z,\infty)}&\le\sup_{z\in\De}\bigg(\frac{\|p\|_L}{|p(z)|}\bigg)^{1/2}=
\bigg(\frac{\|p\|_L}{\dis\min_{z\in\De}|p(z)|}\bigg)^{1/2}\\
&=\bigg(\frac{1+\thi_0}{\Big(\dfrac{\thi_0}{2}\Big)^2-1}\bigg)^{1/2}\Rightarrow\\
\thi_{L,\De}&\le\bigg(\frac{1+\thi_0}{\Big(\dfrac{\thi_0}{2}\Big)^2-1}\bigg)^{1/2}\Rightarrow\
\end{align*}
\[
\rho_L=\thi_L\le\bigg(\frac{1+\thi_0}{\Big(\dfrac{\thi_0}{2}\Big)^2-1}\bigg)^{1/2}.
\]
So, by the above inequality, it suffices to have $\Big(\dfrac{1+\thi_0}{\Big(\dfrac{\thi_0}{2}\Big)^2-1}\Big)^{1/2}<\dfrac{1}{2}$ $(\ast)$, in order to have $\rho_L<\dfrac{1}{2}$.

The minimum natural number $\thi_0$ in order $(\ast)$ holds is $\thi_0=18$. So, for $K_0=\oD(0,1)$, $K_1=\oD(18,1)$, and $\bi_n=2^n$, $n=1,2,\ld$, we have:
\[
\cu(0,\oD,K_1,(2^n))\neq\emptyset.
\]
\end{exm}
\begin{exm}\label{exm3.2}
We fix a positive number $\bi_0\in(0,+\infty)$, $\bi_0\neq1$. We consider the sequence $\bi_n:=\bi^n_0$, $n=1,2,\ld\;$. We fix also a natural number $m_0\ge7$. We give here an example of a compact set $L=\bigcup\limits^{m_0}_{i=0}K_i$, where $K_i$, $i=0,1,2,\ld,m_0$, are pairwise disjoint simple compact sets, $K_0$ is the closed unit disc $z_0=0$, $\vPi=\bigcup\limits^{m_0}_{i=1}K_i$, $K_i$, $i=1,2,\ld,m_0$ are closed discs with the same radius 1, such that $\cu(0,K_0,\vPi,(\bi_0)^n)\neq\emptyset$. $(\ast)$

The centers of the discs $K_i$, $i=1,2,\ld,m_0$ are all on a circle with center 0 and radius $h_0$ for some positive number $h_0>0$ that we will find later. The centers $\rho_j$, $j=0,1,\ld,m_0-1$ of the discs $K_{j+1}$ $j=0,1,\ld,m_0-1$ respectively are vertices of the unique canonical polygon on the circle $C(0,h):=\{z\in\C\mid\,\mid z\mid=h_0\}$ with $m_0$ vertices and one of them is the number $h_0$. Thus, we have
\[
K_{j+1}:=\oD(\rho_j,1), \ \ j=0,1,\ld,m_0-1,
\]
\[
\oD(\rho_j,1):=\{z\in\C\mid\,\mid z-\rho_j\mid\le1\}, \ \ \text{and} \ \ \rho_j=h_0\cdot e^{\frac{2j\pi i}{m_0}}, \ \ j=0,1,\ld,m_0-1,
\]
of course $m_0-1\ge6$.

We suppose that $\bi_0\in(0,1)$. It suffices to have:
\setcounter{equation}{0}
\begin{eqnarray}
\rho_L<\bi_0, \label{eq1}
\end{eqnarray}
because then by Theorem \ref{thm2.2} we take the conclusion in $(\ast)$. So, it suffices to find a compact set $L$ with the above characteristic such that (\ref{eq1}) holds. This means to find the number $h_0$ of course.

We will find the number $h_0=|\rho_j|$ for $j=0,1,\ld,m_0-1$, by imposing sufficient conditions gradually.

First of all we need to have that the discs $K_i$, $i=0,1,\ld,m_0$ are pairwise disjoint. So, if $j_1,j_2\in\{1,2,\ld,m_0\}$, $j_1\neq j_2$, it suffices to have $K_{j_1}\cap K_{j_2}=\emptyset$.

It is enough for this to have
\begin{eqnarray}
h_0>\frac{1}{\sin\Big(\dfrac{\pi}{m_0}\Big)}.  \label{eq2}
\end{eqnarray}
We must have also
\[
K_0\cap K_j, \ \ \text{for every} \ \ j=1,2,\ld,m_0.
\]
It suffices for this to have
\begin{eqnarray}
h_0>2.  \label{eq3}
\end{eqnarray}
The condition $m_0\ge7$ gives that $\dfrac{1}{\sin\Big(\dfrac{\pi}{m_0}\Big)}>2$, so if we take $h_0$ such that (\ref{eq2}) holds then (\ref{eq3}) holds also.

So, if (\ref{eq2}) holds then the discs $K_j$, $j=0,1,\ld,m_0$ are pairwise disjoint.

We fix some positive number $r_0>1$.

We consider the circle
\[
\de_0:=C\bigg(0,\frac{h_0}{2}\bigg):=\bigg\{z\in\C\mid\,\mid z\mid=\frac{h_0}{2}\bigg\}
\]
and the circles $\de_j$, $j=1,2,\ld,m_0$, where
\[
\de_j:=C(\rho_{j-1},r_0):=\{z\in\C\mid\,\mid z-\rho_{j-1}\mid=r_0\} \ \ \text{for} \ \ j=1,2,\ld,m_0.
\]
We have that the circle $\de_j$, $j=0,1,\ld,m_0$ has the compact set $K_j$, $j=0,1,\ld,m_0$ in its interior. We want that everyone from the circles $\de_j$, $j=0,1,\ld,m_0$ has all the others in its exterior. It is easy for this to have
\begin{eqnarray}
r_0<h_0\sin\bigg(\frac{\pi}{m_0}\bigg).  \label{eq4}
\end{eqnarray}
We consider the polynomial
\[
p(z):=z\cdot\prod^{m_0-1}_{j=0}(z-\rho_j).
\]

We take that:
\begin{eqnarray}
\|p\|_L\le(1+h_0)\cdot(1+2h_0)^{m_0-1}.  \label{eq5}
\end{eqnarray}
We also have
\begin{eqnarray}
|p(z)|\ge\bigg(\frac{h_0}{2}\bigg)^{m_0+1} \ \ \text{for} \ \ z\in\de_0 \label{eq6}
\end{eqnarray}
and
\begin{eqnarray}
|p(z)|\ge(h_0-r_0)\cdot r_0\cdot\bigg(2h_0\sin\bigg(\frac{\pi}{m_0}\bigg)-r_0\bigg)^{m_0-1} \ \label{eq7}
\end{eqnarray}
for every $z\in\de_j$, for every $j=1,2,\ld,m_0$.

We take now
\begin{eqnarray}
h_0>\frac{2}{\sin\Big(\dfrac{\pi}{m_0}\Big)}.  \label{eq8}
\end{eqnarray}
We fix
\begin{eqnarray}
r_0:=\frac{1}{2}h_0\sin\bigg(\frac{\pi}{m_0}\bigg).  \label{eq9}
\end{eqnarray}
For such a number $r_0$ as in (\ref{eq9}) we have $1<r_0<h_0\sin\Big(\dfrac{\pi}{m_0}\Big)$.

We set
\begin{eqnarray}
\el_0:=\frac{1}{2}\sin\bigg(\frac{\pi}{m_0}\bigg)\cdot\bigg(1-\frac{1}{2}\sin\bigg(
\frac{\pi}{m_0}\bigg)\bigg)\cdot\bigg(\frac{3}{2}\sin\bigg(\frac{\pi}{m_0}\bigg)\bigg)^{m_0-1}.
\label{eq10}
\end{eqnarray}
Then, we have
\begin{eqnarray}
\bigg(\frac{h_0}{2}\bigg)^{m_0+1}>(h_0-r_0)r_0(2h_0\sin\bigg(\frac{\pi}{m_0}\bigg)-r_0\bigg)^{m_0-1}=
\el_0\cdot h_0^{m_0+1}.  \label{eq11}
\end{eqnarray}
By (\ref{eq6}), (\ref{eq7}) and (\ref{eq11}) we have:
\[
\min_{z\in\De}|p(z)|\ge\el_0\cdot h_0^{m_0+1}, \ \ \text{where}
\]
\begin{eqnarray}
\De=\bigcup^{m_0}_{j=0}\de_j.  \label{eq12}
\end{eqnarray}
By (\ref{eq5}) and (\ref{eq12}) we have:
\begin{eqnarray}
\frac{\|p\|_L}{\dis\min_{z\in\De}|p(z)|}\le
\frac{(1+h_0)(1+2h_0)^{m_0-1}}{\el_0\cdot h_0^{m_0+1}}=
\frac{1}{\el_0}\cdot\frac{1}{h_0}\cdot\bigg(1+\frac{1}{h_0}\bigg)\cdot
\bigg(2+\frac{1}{h_0}\bigg)^{m_0-1}.  \label{eq13}
\end{eqnarray}
Of course
\begin{eqnarray}
\lim_{h\ra+\infty}\bigg(\bigg(1+\frac{1}{h}\bigg)\cdot\bigg(2+\frac{1}{h}\bigg)^{m_0-1}\bigg)=
2^{m_0-1}.  \label{eq14}
\end{eqnarray}
By (\ref{eq13}) and (\ref{eq14}) we have that there exists $h_1>0$ such that
\[
\bigg(1+\frac{1}{h}\bigg)\cdot\bigg(2+\frac{1}{h}\bigg)^{m_0-1}<2^{m_0+1} \ \ \text{and}
\]
\begin{eqnarray}
h>\frac{2}{\sin\Big(\dfrac{\pi}{m_0}\Big)} \ \ \text{for every} \ \ h\ge h_1.  \label{eq15}
\end{eqnarray}
So by (\ref{eq13}) and (\ref{eq15}) we have:
\begin{eqnarray}
\frac{\|p\|_L}{\dis\min_{z\in\De}|p(z)|}<\frac{1}{\el_0}\cdot\frac{1}{h}\cdot2^{m_0+1} \ \ \text{for every} \ \ h\ge h_1.  \label{eq16}
\end{eqnarray}
So by (\ref{eq16}) we take:
\[
\bigg(\frac{\|p\|_L}{\dis\min_{z\in\De}|p(z)|}\bigg)^{\frac{1}{m_0+1}}<
\frac{2}{\sqrt[m_0\!+\!1]{\el_0}}\cdot\frac{1}{h^{1/m_0+1}} \ \ \text{for every} \ \ h\ge h_1.
\]
Because
\[
\frac{2}{\sqrt[m_0\!+\!1]{\el_0}}\cdot\frac{1}{h^{1/m_0+1}}\ra0 \ \ \text{as} \ \ h\ra+\infty,
\]
we can find some positive number $h_2\ge h_1$ such that:
\begin{eqnarray}
\frac{2}{\sqrt[m_0\!+\!1]{\el_0}}\cdot\frac{1}{h^{\frac{1}{m_0+1}}}<\bi_0 \ \ \text{for every} \ \ h\ge h_2.  \label{eq17}
\end{eqnarray}
So, we fix some positive number $h_0\ge h_2$ where $h_2$ is defined as in (\ref{eq17}). Then, we take:
\begin{eqnarray}
\bigg(\frac{\|p\|_L}{\dis\min_{z\in\De}|p(z)|}\bigg)^{\frac{1}{m_0+1}}<\bi_0, \label{eq18}
\end{eqnarray}
where the compact set $L$ is constructed now such that $\rho_j=h_0 e^{\frac{2\pi ji}{m_0}}$ for  $j=0,1,\ld,m_0-1$ from the previous work.

Now let $D:=(\C\cup\{\infty\})\setminus L$. Let $g_D$ be the Green's function for $D$ with pole at $\infty$. By Theorem 5.5.7 (9) of \cite{21} we take:
\begin{eqnarray}
\bigg(\frac{|p(z)|}{\|p\|_L}\bigg)^{\frac{1}{m_0+1}}\le e^{g_D(z,\infty)}, \ \ \text{for every} \ \ z\in D\setminus\{\infty\}. \label{eq19}
\end{eqnarray}
By (\ref{eq18}) and (\ref{eq19}) we take that:
\[
\max_{z\in\De}e^{-g_\De(z,\infty)}<\bi_0\Rightarrow\thi_{L,\De}<\bi_0\Rightarrow
\rho_L<\bi_0
\]
and the desired inequality (\ref{eq1}) holds now. The construction of example is complete now.
\end{exm}
\begin{exm}\label{exm3.3}
In the previous two examples we have given specific examples of sequences $(\bi_n)$ such that their sequences $|\bi_n|^{1/n}$ have a limit point different from 1 and $\cu(\bi_n)\neq\emptyset$.

However, it is desirable to define as many sequences $(\bi_n)$ as we can in this case.

By Theorem \ref{thm2.2} we can succeed this if we define the number $\rho_L$.

We do this here with a specific example.

Let $K_0$ be the hexagon that is defined by coordinates:
\[
-6.5, \ \ -5\pm1.5i, \ \ -5.75\pm2.25i, \ \ -8
\]
and $K_1$ be the square by coordinates $9.5$, $8.75\pm0.75i$, 8 see [9, page 5]. We set $L:=K_0\cup K_1$, $\vPi=K_1=L\sm K_0$. As it is computed in this paper \cite{9} we have $\rho_L:=0.529966\ld$\;.

So by Theorem \ref{thm2.2} we get that for every sequence $(\bi_n)$, such that $|\bi_n|^{1/n}$ has a limit point in $\Big(\rho_L,\dfrac{1}{\rho_L}\Big)$ we have $\cu(\bi_n)=\cu(z,K_0,\vPi,(\bi_n))\neq\emptyset$ for every $z\in\ko_0$.
\end{exm}
\section{The negative case. Cases where we have $\bbb{\cu(\bi_n)=\emptyset}$}\label{sec4}
\noindent

Let some sequence $(\bi_n)$ such that $\bi_n\ra0$. By Theorem \ref{thm2.2} we have that if the sequence $|\bi_n|^{1/n}$ has a limit point in $(\rho_L,1]$ then $\cu(\bi_n)\neq\emptyset$.

So the natural question is what happens when $\underset{n\ra+\infty}{\lim\sup}|\bi_n|^{1/n}\in[0,\rho_L]$.
Consider a compact set $L:=\bigcup\limits^{m_0}_{i=0}K_i$ for some $m_0\in\N$, where $K_i$, $i=0,1,\ld,m_0$ is as in the previous section. We prove here that for such  $L$ there are always cases where $\cu(\bi_n)=\emptyset$. Let $\vPi:=L\sm K_0=\bigcup\limits^{m_0}_{i=1}K_i$ and $z_0\in\ko_0$.
Define $R_0:=dist(z_0,\vPi)$ and $r_0:=dist(z_0,K_0^c)$.
We fix some $w_0\in\vPi$ such that $R_0=|z_0-w_0|$. Of course $\dfrac{r_0}{R_0}\in(0,1)$.

We show the following Proposition \ref{prop4.1}, which implies item $(ii)$ of Theorem \ref{mainresult}.
\begin{prop}\label{prop4.1}
We consider the compact set $L$, and the number $w_0$ as previously and some sequence $(\bi_n)$ of complex numbers such that $\underset{n\ra+\infty}{\lim\sup}|\bi_n|^{1/n}\in\Big[0,\dfrac{r_0}{R_0}\Big)$. Then we have
\[
\cu(\bi_n)=\cu(z_0,K_0,\vPi,(\bi_n))=\emptyset.
\]
\end{prop}
\begin{Proof}
It is obvious that $\bi_n\ra0$. Suppose that $\bi_n\neq0$ for every $n=1,2,\ld\;$
The other case is similar, so it is omitted .
We set $\Si:=\underset{n\ra+\infty}{\lim\sup}|\bi_n|^{1/n}\in\Big[0,\dfrac{r_0}{R_0}\Big)$ and
assume that $\cu(z_0,K_0,\vPi,(\bi_n))\neq\emptyset$. Let $f\in\cu(z_0,K_0,\vPi,(\bi_n))$
and $\el_0$ be the radius of convergence of the Taylor's development of $f$ with center $z_0$. Of course, $\el_0\ge r_0$. We choose some $\e_0\in(0,r_0)$ such that: $\Si\in\Big[0,\dfrac{r_0-\e_0}{R_0}\Big)$ (for example take $\e_0\in(0,r_0-\Si R_0)$.
Let $N_0\in\N$. Then, 
\setcounter{equation}{0}
\begin{eqnarray}
|\bi_{N_0}S_{N_0}(f,z_0)(w_0)| \le|\bi_{N_0}|\bigg(\frac{R_0}{r_0-\e_0}\bigg)^{N_0}\cdot\sum^{N_0}_{k=0}
\bigg|\frac{f^{(k)}(z_0)}{k!}\bigg|(r_0-\e_0)^k.  \label{eq1}
\end{eqnarray}
Because $r_0-\e_0\in(0,\el_0))$ we have:
\begin{eqnarray}
\sum^{+\infty}_{k=0}\bigg|\frac{f^{(k)}(z_0)}{k!}\bigg|(r_0-\e_0)^k=A\in[0,+\infty). \label{eq2}
\end{eqnarray}
By (\ref{eq1}) and (\ref{eq2}) we get:
\begin{eqnarray}
|\bi_{N_0}\cdot S_{N_0}(f,z_0)(w_0)|\le|\bi_{N_0}|\cdot\bigg(\frac{R_0}{r_0-\e_0}\bigg)^{N_0}\cdot A \ \ \text{for} \ \ N_0\in\N. \label{eq3}
\end{eqnarray}
Consider some positive number
\begin{eqnarray}
\thi_0\in\bigg(\Si,\frac{r_0-\e_0}{R_0}\bigg).  \label{eq4}
\end{eqnarray}
Then, there exists some natural number $n_0$ such that
\[
|\bi_n|^{1/n}<\thi_0 \ \ \text{for every} \ \ n\ge n_0.
\]
So
\begin{eqnarray}
|\bi_n|<\thi^n_0 \ \ \text{for} \ \ n\ge n_0.  \label{eq5}
\end{eqnarray}
By (\ref{eq3}), (\ref{eq4}) and (\ref{eq5}) we have:
\begin{eqnarray}
|\bi_n\cdot S_n(f,z_0)(w_0)|<\bigg(\frac{\thi_0R_0}{r_0-\e_0}\bigg)^n\cdot A \ \ \text{for every} \ \ n\ge n_0 \ \ \text{where} \ \ \frac{\thi_0R_0}{r_0-\e_0}\in(0,1)  \label{eq6}
\end{eqnarray}
 and (\ref{sec36}) implies: $\dis\lim_{N\ra+\infty}\bi_NS_N(f,z_0)(w_0)=0$ that means that $f\notin\cu(z_0,K_0,\vPi,(\bi_n))$. This is a contradiction, so we have proved Proposition \ref{prop4.1}. \qb
\end{Proof}

We may have $\cu(\bi_n)=\emptyset$ when the sequence $\bi_n$ tends to infinity very fast.
Let $L=\bigcup\limits^{m_0}_{i=0}K_i$ some compact set.
We fix some $z_0\in\ko_0$ and we set $r_0:=dist(z_0,K^c_0)$, $\vPi:=L\sm K_0=\bigcup\limits^{m_0}_{i=1}K_i$ and $\OO_1:=(\C\cup\{\infty\})\sm\vPi$. Let $g_{\OO_1}$ be the Green's function for $\OO_1$ with pole at infinity. Of course,
\[
D_1:=\oD(z_0,r_0)=\{z\in\C\mid |z-z_0|\le r_0\}\subseteq K_0\subset\OO_1\sm\{\infty\}.
\]
We set $M_0:=e^{\dis\max_{z\in D_1}g_{\OO_1}(z,\infty)}$.

Using the above notations we prove the following proposition, which implies item $(iii)$ of Theorem \ref{mainresult}.
\begin{prop}\label{prop4.2}
Let $(\bi_n)$ be a sequence of complex numbers such that\\ $\underset{n\ra+\infty}{\lim\inf}|\bi_n|^{1/n}>M_0$. Then $\cu(z_0,K_0,\vPi,(\bi_n))=\emptyset$.
\end{prop}
\begin{Proof}
Of course $M_0>1$ and so $\bi_n\ra\infty$. Without loss of generality we can suppose that $\bi_n\neq0$, for every $n=1,2,\ld\;.$
We suppose that
\[
\cu(\bi_n)=\cu(z_0,K_0,\vPi,(\bi_n))\neq\emptyset.
\]
Let $f\in\cu(\bi_n)$. Of course $f\neq\O$, where $\O:K_0\ra\C$, $\O(z)=0$, $\fa\,z\in K_0$.
So there exists $z_1\in D(z_0,r_0)$ such that $f(z_1)\neq0$ or else by the principle of identity we would have $f=\O$ that is false. By the continuity of $f$ we have that there exist $\thi_1>0$ and $\de_1>0$ such that $D_{\de_1}:=\overline{D}(z_1,\de_1)\subset D(z_0,r_0)$ and $|f(z)|>\thi_1$ for every
\setcounter{equation}{0}
\begin{eqnarray}
z\in D_{\de_1}.  \label{eq1}
\end{eqnarray}
Of course we have:
\begin{eqnarray}
\|S_n(f,z_0)-f\|_{D_{\de_0}}\ra0 \ \ \text{as} \ \ n\ra+\infty.  \label{eq2}
\end{eqnarray}
Let some $\thi_2\in(0,\thi_1)$. By (\ref{eq1}) and (\ref{eq2}) there exists some $n_1\in\N$, such that:
\begin{eqnarray}
|S_n(f,z_0)(z)|>\thi_2 \ \ \text{for every} \ \ z\in D_{\de_0}, \ \ n\ge n_1.  \label{eq3}
\end{eqnarray}
We denote $m_n:=deg S_n(f,z_0)$, $n\ge n_1$. Of course $m_n\le n$, $\fa\,n\ge n_1$.

By Bernstein's Lemma (Theorem 5.5.7 (9) \cite{21}) for the polynomial $S_n(f,z_0)$, $n\ge n_1$ we get
\begin{eqnarray}
\bigg(\frac{|S_n(f,z_0)(z)}{\|S_n(f,z_0)\|_\vPi}\bigg)^{1/m_n}\le
e^{g_{\OO_1}(z,\infty)} \ \ \text{for} \ \ z\in D_{\de_0}.  \label{eq4}
\end{eqnarray}
By (\ref{eq4}) and the definition of the number $M_0$ we get:
\begin{eqnarray}
|\bi_nS_n(f,z_0)(z)|^{1/n}\le M_0^{\frac{m_n}{n}}\Big\|\bi_nS_n(f,z_0)\Big\|_\vPi^{1/n}, \ \ n\ge n_1, \ \ z\in D_{\de_0}.  \label{eq5}
\end{eqnarray}
By the supposition of our proposition it follows that:
\[
\underset{n\ra+\infty}{\lim\inf}^{1/n}>M_0.
\]
Choose some $M_1\in\big(M_0,\underset{n\ra+\infty}{\lim\inf}|\bi_n|^{1/n}\big)$.
There exists some $n_2\ge n_1$ such that:
\begin{eqnarray}
|\bi_n|^{1/n}>M_1 \ \ \text{for} \ \ n\ge n_2.  \label{eq6}
\end{eqnarray}
By (\ref{eq3}), (\ref{eq5}) and (\ref{eq6}) we have
\begin{eqnarray}
M_1\thi_2^{1/n}<M^{\frac{m_0}{n}}_0\Big\|\bi_nS_n(f,z_0)\Big\|_\vPi^{1/n}, \ \ n\ge n_2  \label{eq7}
\end{eqnarray}
and (\ref{eq7}) implies
\begin{eqnarray}
\frac{M_1}{M_0}\thi^{1/n}_2<\Big\|\bi_nS_n(f,z_0)\Big\|^{1/n}_\vPi, \ \ n\ge n_2.  \label{eq8}
\end{eqnarray}
Because $M_1>M_0$, by (\ref{eq8}) there exists $n_3\ge n_2$ such that:
\[
\|\bi_nS_n(f,z_0)\|_\vPi>1, \ \ n\ge n_3.
\]
This means that $f\notin\cu(\bi_n)$, that is false, and the proof is complete. \qb
\end{Proof}
\begin{rem}\label{rem4.3}
Theorem \ref{thm2.2} of this paper paper and Propositions 3.1 and 2.3 of \cite{24} give an alternative proof of the main result in \cite{24}.
\end{rem}
\vspace*{1cm}
Tsirivas Nikos \\
Postdoctoral researcher \\
University of Crete \\
Department of Mathematics and Applied Mathematics\\
Panepistimiopolis Voutes, 700-13, Heraklion, Crete, Greece\\
email: tsirivas@uoc.gr

\end{document}